\documentclass[regno, reqno,11pt]{amsart}
\usepackage{amsmath,amsthm,amssymb,amsfonts,epsfig,color,graphicx,enumerate,accents,subfigure,comment, esint}
\usepackage[a4paper,margin=2.59truecm]{geometry}
\usepackage{graphicx}
\usepackage{enumitem}

\usepackage[linkcolor=blue,colorlinks=true]{hyperref}\usepackage{fancyhdr}
\DeclareMathOperator{\dive}{div}

\DeclareMathOperator{\Tr}{Tr}

\def\eps{{\varepsilon}}

\def\N{\mathbb{N}}

\def\R{\mathbb{R}}
\def\O{\Omega}
\def\o{\omega}

\def\HH{\mathcal{H}}

\def\DD{\mathcal{D}}

\def\FF{\mathcal{F}}

\def \< {\langle}
\def \>{\rangle}
\def\HH{\mathcal{H}}
\def\II{\mathcal{I}}

\def\KK{\mathcal{K}}

\def\MM{\mathcal{M}}
\def\PP{\mathcal{P}}

\def\TT{\mathcal{T}}

\def\f{\phi}
\def\k{\kappa}

\newcommand{\be}{\begin{equation}}
\newcommand{\ee}{\end{equation}}

\newcommand{\bs}{\begin{split}}
\newcommand{\es}{\end{split}}

\newcommand{\Linfty}[2]{\| #1 \|_{L^{\infty}(#2)} }
\renewcommand{\L}[3]{\| #2 \|_{L^{#1}(#3)}}
\renewcommand{\b}[1]{\bar{#1}}

\numberwithin{equation}{section}
\theoremstyle{plain}
\newtheorem{theorem}{Theorem}[section]
\newtheorem{lemma}[theorem]{Lemma}
\newtheorem{corollary}[theorem]{Corollary}
\newtheorem{proposition}[theorem]{Proposition}
\newtheorem{definition}[theorem]{Definition}
\theoremstyle{remark}
\newtheorem{remark}[theorem]{Remark}

\title[Transmission problem tangential touch]{ Tangential contact of free boundaries and the fixed boundary for variational solutions to a free transmission Problem}

\author{Diego Moreira$^{*}$}
\address{$^*$Departamento de Matemática, Universidade Federal da Ceara(Fortaleza, Brazil)}
\email{$^*$dmoreira@mat.ufc.br}
\author{Harish Shrivastava$^\dagger$}
\address{$^\dagger$Tata Institute of Fundamental Researcher-Centre of of Applicable Mathematics}
\email{$^\dagger$harish21@tifrbng.res.in, }

\begin{document}

\begin{abstract}
In this article we study functionals of the following type 
$$
\int_{\O} \Big ( \< A(x,u)\nabla u, \nabla u\> + \Lambda (x,u) \Big )\,dx
$$
here $A(x,u)= A_+(x)\chi_{\{u>0\}}+A_-(x) \chi _{\{u\le0\}}$ for some elliptic and bounded matrices $A_{\pm}$ with Hölder continuous entries and $\Lambda(x,u) = \lambda_+(x) \chi_{\{u>0\}} + \lambda_-(x) \chi_{\{u\le 0\}}$.{ We prove that the free boundaries of minimizers of the above functional touches the fixed boundary $\partial \O$ in a tangential fashion, provide the graph of boundary data touches its zeros smoothly. This assumption is reflected in the \eqref{DPT} condition.}
\end{abstract}
\medskip 

\maketitle


\textbf{2010 Mathematics Subjects Classification:} 49J05, 35J20, 35A15, 35B65, 35R35.

\tableofcontents

\section{Introduction}\label{introduction}
Phenomenas which involve interaction of one medium with another can be modeled through free transmission problems of the form \eqref{functional} and \eqref{PDE formulation}. This class of problems appear naturally when each medium in the mixture follows its own diffusion laws. Few examples of phenomenas that can be modeled through transmission problems are mixture of different chemicals, conductivity (electric and thermal) of composite materials or a material operating close to thresholds like triple point, melting point or breakdown potential. It is indeed a challenging task to model the properties of certain materials, particularly when they are close to threshold points which lead to abrupt changes in their diffusion properties. Applications and results dealing with such models can be found in \cite{M10}, \cite{AI96} and references therein.

As nicely pointed out in \cite{CSS21},  M. Picone introduced transmission problems in theory of elasticity in 1954 (\cite{picone}) and the theory was further developed by Lions \cite{L56}, Stampaccia \cite{S56} and Campanato \cite{C57}. The non-divergence case was treated by Schechter in 1960 (\cite{S60}). For further details, we refer to \cite{CSS21}, from where we learn about the history of those developments.  
 
Roughly speaking, the models discussed above can be represented by the PDE
\be\label{PDE formulation}
\dive A(x,u,\nabla u) = f(x,u)
\ee
where $A(x,u,\nabla u)$ and $f(x,u)$ have jump discontinuities with respect to the variable $u$. 
Transmission problems can also be posed in variational setting where configurations corresponding to least energy are of interest.  We study minimizers of functionals of the following type
\be\label{variational formulation}
\int_{\O}H(x,u,\nabla u)\,dx
\ee
where $H(x,u,\nabla u)$ has a jump discontinuity with respect to the variable $u$. Some developments in above mentioned variational formulation \eqref{variational formulation} can be found in \cite{HS21}, \cite{MS21} and \cite{TA15}. Moreover, the form \eqref{PDE formulation} is also studied in \cite{KLS17}, \cite{KLS19}. 

{Several authors have extended the works of Alt, Caffarelli, Friedman \cite{ACF84} on free boundary problems with constant coefficients to the case of variable coefficients. The works of Argiolas, Ferrari, Sandro, Cerutti \cite{F06}, \cite{AF09}, \cite{FS07}, \cite{CFS04} put an effort to extend the seminal works of Caffarelli \cite{C1}, \cite{C2} to the case of variable coefficients. The fundamental paper by Caffarelli \cite{C3} deals with free boundary regularity for fully non-linear PDEs with variable coefficients. The series of papers by De Silva, Salsa, Ferrari \cite{DFS19}, \cite{DFS16}, \cite{DFS15.2}, \cite{DFS15.1}, \cite{DFS14.1}, \cite{DFS14} extends the results due to Caffarelli \cite{C1}, \cite{C2}, \cite{C3} to non-homogeneous case. One can refer to work of Ferrari, Salsa \cite{FS07} for free boundary regularity for non-divergence PDE with variable coefficients and drift. G. David, Engelstein, Garcia, Toro \cite{DEGT19} look into almost minimizers of Bernoulli type free boundary problem with variable coefficient. See also the recent results by Ferrari, Lederman \cite{FL21} for the case of variable exponents. }

{All the above mentioned works are the interior regularity results on solutions and their free boundaries. To the best of our knowledge, there are very few known results which deal with up to the boundary regularity in FBPs (for example the papers \cite{KKS06}, \cite{KS15}, \cite{MS22}, \cite{GL19} deal with contact geometry of free boundary and fixed boundar and \cite{CS19} looks into upto the boundary regularity of solutions). We refer to our previous work \cite{MS22}, which proves the tangential interaction of free boundaries and fixed boundary in the case of variable coefficient Bernoulli type free boundary problem.}

{In this paper, we explore a slightly more delicate aspect of Bernoulli type free boundary problems: the case of variable coefficients in setting of free transmission problem \eqref{variational formulation} (i.e. coefficients with a jump discontinuity along the free boundaries).  
We try to add to the understanding of solutions and their free boundaries in this setting. We look at the interaction of the free boundaries of solutions and the fixed boundary of domain for functionals of the form \eqref{functional}.}

Main result proven in this work is the free boundaries of minimizers of functional of the form \eqref{functional} touches the fixed boundary in a tangential fashion, provided the boundary data is well behaved (c.f. \eqref{DPT} condition). Our findings are in line with \cite{KKS06}, \cite{MS22} and \cite{KS15}. In \cite{KS15}, authors show that the interaction of free boundary and fixed boundary is transversal in absence of boundary data with \eqref{DPT}.

As discussed above, interactions of respective diffusions of composite materials can be modeled by free transmission problems.  In certain circumstances, it is interesting to study the free boundary of solution and fixed boundary in their contact set (if exists). The Dam problem \cite{AG82} and Jets, Wakes and Cavities \cite{ZG57} model phenomenas.

Recent works of Indrei \cite{I00}, \cite{I19} study such interactions for fully non-linear obstacle problems. We also refer to \cite{GL19} which sheds more light into angle of contact between fixed boundary and free boundary for one phase Bernoulli problem. Since we are dealing with a free transmission problem, we make use of a very handy tool, particularly useful in study of such scenarios: the $\TT_{a,b}$ operator. It was introduced in \cite{MS21}. The $\TT_{a,b}$ operator bridges between free transmission problems of the form \eqref{J0} and functionals for Alt-Caffarelli-Friedman type \cite{ACF84}. A detailed discussion on $\TT_{a,b}$ operator can be found in \cite[Section 3]{MS21}.

We follow the strategy of classifying global profiles. Our paper extends the result of \cite{KKS06} to the case of free transmission problems with Hölder continuous coefficients. In order to study the blow-ups, we establish non-degeneracy and energy estimates for solutions to transmission problems which may be of independent interest to the readers. {The fact that we deal with free transmission problems, i.e. problem with discontinuous coefficients, brings new difficulties in to our setting. For instance the technique of using harmonic replacements (used in \cite[Lemma 3.7]{MS22}, \cite[Theorem 3.1]{KKS06}) does not appear to be available in our setting. Therefore, we bring a new, more suitable approach to prove the energy estimates (c.f. Theorem \ref{Energy Estimates}) which was obtained via harmonic replacement in \cite[Lemma 3.7]{MS22}, \cite[Theorem 3.1]{KKS06} in their corresponding settings. We employ Widman's hole filling technique to prove the uniform energy estimates. Our strategy in the proof of  Theorem \ref{Energy Estimates} may be of independent interest to the readers.}

{Moreover, in the proof of non-degeneracy near the free boundary (c.f. Proposition \ref{non degeneracy}), we construct delicate variational barriers. In order to control these barriers, we use Gruter-Widman estimates on Green's functions \cite{GW82} along the positive phase of our variational solution. The compactness lemma (Proposition \ref{blowup convergence}) is proven by careful choice of  a test function (c.f. \eqref{test function}) which allows u to pass through the limit while preserving the minimality condition derived from nature of blowups of a minimizer.}

In Section \ref{setting the problem}, we introduce the problem and present definitions and notations to be used in the rest of the paper. We also present the main Theorem \ref{tangential touch}. In Section \ref{section 3}, we prove compactness, energy estimates, non-degeneracy and other supporting lemmas for blowups of solutions to \eqref{functional}.  At last in Section \ref{main proof}, we prove our main result Theorem \ref{tangential touch}.

\section{Setting up the problem}\label{setting the problem}

Objective of the paper of to study the behaviour of free boundary near the boundary of domain, for minimizers Bernouli type free transmission problems with Hölder continuous matrix coefficients. 
\be\label{functional}
J(v;A_{\pm},\lambda_{\pm},Q,\O):=\int_{\O} \Big ( \< A(x,u)\nabla v,\nabla v \>+ Q(x) \Lambda(v)\Big )\,dx.
\ee
Here, $A(x,s) := A_+(x) \chi_{\{s>0\}} + A_-(x) \chi_{\{s\le 0\}}$ and $\Lambda(s):= \lambda_+\chi_{\{s>0\}}+ \lambda_{-}\chi_{\{s\le 0\}}$. In order to study the contact of free boundaries and fixed boundaries, we consider the minimization problem in the domain $B_R^+$ denoted as

\[
\begin{split}
B_R^+ := \Big \{ x\in B_R\mbox{ such that } x_N>0 \Big \}\\
B_R' := \Big \{ x\in B_R\mbox{ such that } x_N=0 \Big \}.
\end{split}
\]
For $x\in \R^N$ we denote $x'\in \R^{N-1}$ as the projection of $x$ on the plane $\{x_N=0\}$, we  denote the tangential gradient of a function $u$, $\nabla' u$ as follows
$$
\nabla' u := \Big ( \frac{\partial u}{\partial x_1}, ...\,\frac{\partial u}{\partial x_{N-1}} \Big ).
$$
For $\f\in H^1(B_R^+)$, we define the affine space $H^1_\f(B_R^+)$ as follows,
\be\label{Kf}
H^1_\f(B_R^+) = \left \{ v\in H^1(B_R^+) \,:\, v-\f\in H_0^1(B_R^+) \right \}.
\ee
For a given function $v\in H^1(B_R^+)$, we denote the set $F(v)$ as $F(v):=\partial \{v>0\}\cap B_R^+$ and $Id$ is notation for $N\times N$ identity matrix. Throughout the paper, the dependence of a function $f(\cdot,u)$ on the variable $u$ will be seen as follows 
$$
f(\cdot,u):= f(\cdot,u^+)\chi_{\{u>0\}} + f(\cdot,u^-)\chi_{\{u\le 0\}}.
$$

\begin{definition}\label{with density}
A function $u\in H^1(B_{2/r}^+)$ is said to belongs to the class $\MM_{r.\DD}(\alpha, M,\lambda_{\pm},\mu,q,\o)$ if there exists $A_{\pm}\in C^{\alpha}(B_{2/r}^+)^{N\times N}$, $\f \in C^{1,\alpha }(B_{2/r}^+)$, $0<\mu,q<1$, $\DD>0$ and a modulus of continuity $\o:\R^+ \to \R^+$ such that 
\begin{enumerate}[label=\textbf{$($M\arabic*$)$}]
\item \label{P1}$\Linfty{A_{\pm}}{B_{2/r}^+},\,\Linfty{ \nabla \f}{B_{2/r}^+}\le M$, $ [A_{\pm}]_{C^{\alpha}(B_{2/r}^+)},[ \nabla \f] _{C^{0,\alpha}(B_{2/r}^+)}\le r^{\alpha}M$ and \newline $|\f(x')| \le Mr^{1+\alpha}|x'|^{1+\alpha}$ $($$x'\in B_{2/r}'$$)$. $\f$ satisfies the following Degenerate Phase Transition condition \eqref{DPT} mentioned below.
\be\label{DPT}\tag{DPT}
\mbox{$\forall x'\in B_{2/r}'$ if $\f (x')=0$, then $|\nabla' \f(x')|=0$}.
\ee
\item \label{P2} $\mu |\xi|^2 \le \< A_{\pm}(x)\xi,\xi\> \le \frac{1}{\mu}|\xi|^2$ for all $x\in B_{2/r}^+$ and $\xi \in \R^N$. There exist $\mu \le a_{\pm} \le \frac{1}{\mu}$ such that $A_{\pm}(0)=a_{\pm} Id$.
\item \label{P4} $u$ minimizes $J(\cdot; A_{\pm},\lambda_{\pm},Q,B_{2/r}^+)$ $($c.f. \eqref{functional}$)$ that is for every $u-v\in H_0^1(B_{2/r}^+)$
$$
\int_{B_{2/r}^+} \Big ( \< A(x,u) \nabla u,\nabla u  \>+ Q(x)\Lambda (u) \Big ) \,dx \le \int_{ B_{2/r}^+} \Big ( \< A(x,u) \nabla v,\nabla v  \>+ Q(x)\Lambda (v) \Big )\,dx \;\;
$$ and $0\in F(u)$.
\item \label{P3} $0<\lambda_-<\lambda_+$ and $Q\in C(\overline {B_{2/r}^+})$ with modulus of continuity $\o$, \newline i.e. $|Q(x)-Q(y)| \le \o(r|x-y|)$ for $x,y\in B_{2/r}^+$. Moreover $0<q\le Q(x) \le \frac{1}{q}$ for all $x\in B_{2/r}^+$.
\item \label{P5} $u \in H_{\f}^1(B_{2/r}^+)$. 
\item \label{P6} There exists $r_0>0$ such that for all $0<\rho\le r_0$ we have 
\be\label{density}
\frac{|B_\rho^+(0)\cap \{u>0\}|}{|B_\rho^+(0)|}>\DD.
\ee
 
\end{enumerate}
\end{definition}

\begin{remark}\label{remark on DPT}
If $\f\in C^{1,\alpha}(B_2^+)$ satisfies (\ref{DPT}), from \cite[Lemma 10.1]{BM21}, we know that $\f^{\pm} \in C^{1,\alpha}(B_2')$ and also
$$
\| \f^{\pm} \|_{C^{1,\alpha}(B_2')} \le \| \f \|_{C^{1,\alpha}(B_2')}.
$$
\end{remark}
Given $v\in H^1(B_R^+)$ and $r>0$, we define the blow-up $v_r\in H^1(B_{R/r}^+)$ as follows  
\be\label{blowup}
v_r(x):=\frac{1}{r}v(rx).
\ee
For the coefficient matrix $A$, $A^r(x)$ is defined as follows
\be\label{Ar}
A^r(x):=A(rx).
\ee
and similarly we define 
\be\label{Qr}
Q^r(x) = Q(rx)
\ee
\begin{remark}\label{rescaling invariance}
One can check that if $u\in \MM_{1}$, then $u_r \in \MM_{r}$. Indeed if $u\in \MM_{1}$ and $u$ minimizes the functional $J$ (c.f. \ref{P4})
$$
J(v;A_{\pm},\lambda_{\pm},Q,B_2^+):= \int_{B_2^+}\Big ( \< A(x,v)\nabla v,\nabla v \> + Q(x)\Lambda(v)\Big )\,dx,
$$
with boundary data $\f\in C^{1,\alpha}(B_2^+)$ (i.e. $u\in H_{\f}^1(B_2^+)$). Then by simple change of variables we can check that $u_r\in H^1_{\f_r}(B_{2/r}^+)$ (this verifies \ref{P5}) and $u_r$ minimizes 
$$
J(v;A_{\pm}^r,\lambda_{\pm},Q^r,B_{2/r}^+):= \int_{B_{2/r}^+}\Big ( \< A^r(x)\nabla v,\nabla v \> + Q^r(x)\Lambda(v)\Big )\,dx,\qquad  (\Lambda(s)=\lambda_+ \chi_{\{s>0\}}+\lambda_- \chi_{\{s\le0\}}  ).
$$
Moreover, if $A$ and $\f$ satisfy the conditions \ref{P1}, \ref{P2} for $r=1$, then $A^r$, $Q^r$ and $\f_r$ satisfy \ref{P1}, \ref{P2} for $r$. \ref{P3} and \ref{P6} remain invariant under the change variables. Therefore $u_r\in \MM_r$. 
\end{remark}
In order to study the blow-up limits ($\lim_{r\to 0} u_r$) of functions $u\in \MM_1$, we define a class of global solutions $\MM_{\infty}(C,\lambda_{\pm},\mu)$. Let us set the following notation before
$$
\Pi= \{x\,:\, x_N=0\}.
$$

\begin{definition}[Global solution]\label{global solution}
A function $u\in H^1(\R^N_+)$ belongs to the class $\MM_{\infty}(C, \lambda_{\pm}, \mu)$, that is, $u$ is a global solution if there exists $C>0$, $\mu\le a_{\pm}\le \frac{1}{\mu}$ and $0<\lambda_+<\lambda_-$ such that
\begin{enumerate}[label=\textbf{$($G\arabic*$)$}]
\item \label{G1} $|u(x)|\le C|x|$, for all $x\in \R^N_+$,
\item \label{G2}$u$ is continuous up to the boundary $\Pi$,
\item \label{G3} $u=0$ on $\Pi$,
\item \label{G4} and for every ball $B_r(x_0)$, $u$ is a minimizer of $J(\cdot; {a_{\pm},\lambda_{\pm},B_r(x_0)\cap \R^N_+})$ $($c.f. \eqref{functional}$)$, that is 
$$
\int_{B_r(x_0)\cap \R^N_+} \Big (a(u) |\nabla u|^2 +  \Lambda(u)\Big )\,dx \le \int_{B_r(x_0)\cap \R^N_+} \Big (a(u) |\nabla v|^2 +  \Lambda(v) \Big )\,dx
$$
for every $v\in H^1(B_r(x_0)\cap \R^N_+)$ such that $u-v\in H_0^1(B_r(x_0)\cap \R^N_+)$. \newline Here  $a(u) := a_+ \chi_{\{u>0\}} + a_- \chi_{\{u\le 0\}}$ and $\Lambda(u) := \lambda_+ \chi_{\{u>0\}} + \lambda_- \chi_{\{u\le 0\}}$. We denote the class $\MM_{\infty}(C,\lambda_{\pm},1)$ as $\PP_{\infty}(C,\lambda_{\pm})$ $($c.f. \cite[Definition 2.3]{KKS06}, \cite[Definition 2.3]{MS22}$)$.
\end{enumerate}
\end{definition}

\begin{remark}
In the absence of ambiguity on values of $\alpha, M,\lambda_{\pm}, \DD,\mu,q,\o$ we use the notation  $\MM_{r}$ in place of $\MM_{r,\DD}(\alpha, M,\lambda_{\pm}, \mu,q,\o)$. Similarly, we denote $\MM_{\infty}(C,\mu,\lambda_{\pm})$ as $\MM_{\infty}$ when the parameters $C,\mu,\lambda_{\pm}$ are unambiguous.
\end{remark}
\begin{definition}\label{Tab}
Let $a,b>0$ and $D$ be an open Lipschitz set, we define $\TT_{a,b}:H^{1}(D)\to H^{1}(D)$ as follows
$$
\TT_{a,b}(v) := a v^+\,-\,bv^-.
$$
\end{definition}

\begin{remark}\label{prop of Tab}
Here we list some essential properties of the operator $\TT_{a,b}$. One can refer to \cite[Section 3]{MS21} for a detailed exposition. For $v\in H^1(\O)$ the following holds, 
\begin{itemize}
\item  $v\ge 0$ if and only if $\TT_{a,b}(v) \ge 0$. More precisely,
\be\label{unique pm decomp}
\begin{split}
[\TT_{a,b}(v)]^+& =  av^+=(av)^+\mbox{ a.e. in $\O$}\\
[\TT_{a,b}(v)]^- &= b v^-= (bv)^-\mbox{ a.e. in $\O$}.
\end{split}
\ee
\item We define $\lambda (v) := \lambda_{+}\chi_{\{v>0\} \cap \O} + \lambda_- \chi_{\{v\le 0\}\cap \O}$, then $\Lambda(\cdot)$ remains invariant under $\TT_{a,b}$ operator.
\be\label{lambda preserved}
\Lambda(v) = \Lambda(\TT_{a,b}(v)) \;\mbox{a.e. in $\O$}.
\ee
\end{itemize}
\end{remark}

\subsection{The main Result}

The main result of this paper is the following.

\begin{theorem}[Tangential touch for minimizers] \label{tangential touch}
Let $u\in \MM_1$. Then there exists a modulus of continuity $\sigma$ $($$\sigma$ non-decreasing, continuous and $\sigma(0)=0$$)$ and a constant $\rho_0>0$ such that
$$
F(u)\cap B_{\rho_0}^+ \subset \{ x: \; x_N\le \sigma(|x|)|x| \}.
$$  
The function $\sigma$ and constant $\rho_0$ are independent of the choice of $u$ and depends only on the parameters $\alpha, M,\lambda_{\pm},\mu$ and $\DD$.
\end{theorem}

\section{The tangential touch of free boundary with fixed boundary}\label{section 3}
\subsection{Boundary regularity.}
The following lemma is a classical result, c.f. \cite[Remark 4.2]{AC81} and \cite[Lemma 3.1]{MS22}.  

\begin{lemma}$($\cite[Remark 4.2]{AC81}, \cite[Lemma 3.1]{MS22}$)$\label{subharmonic}
Any non-negative continuous function \newline $w \in C(\O)$ such that $\dive(A(x)\nabla w)=0$ weakly in $\{w>0\}\cap \O$ for some strictly elliptic and bounded matrix $A$, then $w\in H^1_{loc}(\O)$ and $\dive(A(x)\nabla w)$ is a non-negative in weak sense.

\end{lemma}
\begin{lemma}\label{Tab blowup are global}
Let $u_0\in \MM_{\infty}$ if and only if $\TT_{\sqrt{a_+},\sqrt{a_-}}(u_0)\in \PP_{\infty}$.
\end{lemma}
\begin{proof}
Lemma \ref{Tab blowup are global} is a direct consequence of \cite[Lemma 3.10]{MS21}.
\end{proof}

We can prove boundedness and Hölder continuity of minimizers by well known results developed in \cite{eg05}.

\begin{lemma}\label{holder}
For every $u\in \MM_r$ there exists $0<\alpha_0<1$ such that $u\in C^{\alpha_0}(\overline {B_2^+})$. Here $\alpha_0:= \alpha_0( M,\lambda_{\pm},\mu,q)$ and we have the estimates 
$$
\|u\|_{C^{\alpha_0}(B_2^+)} \le C(\mu, q, \lambda_{\pm})M.
$$
\end{lemma}
\begin{proof}
We note that the functional $J(\cdot;A_{\pm},\lambda_{\pm},Q,B_2^+)$ satisfies the assumptions of \cite[Theorem 7.3]{eg05} and the boundary data $\phi\in C^{1,\alpha}(\partial {B_2^+})$. Thus, all the assumptions mentioned the discussion in \cite[Section 7.8]{eg05} are satisfied, therefore from the discussions therein, we have  $u\in C^{\alpha_0}(\overline {B_2^+})$ for some $0<\alpha_0 <1$ and 
$$
\|u\|_{C^{\alpha_0}(\overline{B_2^+})} \le C(\mu,q,\lambda_{\pm})\Linfty{u}{B_2^+}.
$$
From Lemma \ref{subharmonic} and comparison principle, we have $\Linfty{u^{\pm}}{B_2^+} \le \Linfty{\f^{\pm}}{\partial B_2^+}\le M$, therefore $\Linfty{u}{B_2^+}\le M$. This proves Lemma \ref{holder}.
\end{proof}

As a consequence of continuity of $u\in \MM_1$ from Lemma \ref{holder}, $\{u>0\}\cap B_2^+$ and $\{u<0\}\cap B_2^+$ are open sets. Since $u^{\pm}$ are $A_{\pm}$-harmonic in $\{u^{\pm}>0\}\cap B_2^+$, we have the following corollary from Lemma \ref{subharmonic}.
\begin{corollary}[Corollary to Lemma \ref{subharmonic} and Lemma \ref{holder}]\label{A subharmonic}
$u^{\pm}$ are $A_{\pm}$-subharmonic in $B_2^+$. 
\end{corollary}

\begin{lemma}[Linear growth]\label{l}
If $u \in \MM_1$. Then we have 
\be\label{lg}
|u(x)| \le C(\mu, [A_{\pm}]_{C^{\alpha}(B_2^+)},N) M |x|\qquad \forall x\in B_1^+.
\ee
\end{lemma}

\begin{proof}
Lemma \ref{l} is proven in \cite[Lemma 3.5]{MS22}, we provide a sketch in this paper. 

We consider the harmonic replacement for $u^+$ in $B_2^+$. Let $w$ be such that 
$$
\begin{cases}
\dive(A_+(x)w) = 0,\qquad \mbox{in $B_2^+$}\\
w = \f^+,\qquad \qquad\qquad\mbox{in $\partial B_2^+$}.
\end{cases}
$$
Since $u^+$ is $A_+$-subharmonic (c.f. Lemma \ref{subharmonic}). By comparison principle, 
\be\label{borsuk 2}
\begin{split}
u^+(x)\le w(x)&\le (\|\nabla w\|_{L^{\infty}(B_1^+)} + M)|x|,\qquad x\in B_1^+.
\end{split}
\ee
From \cite[Theorem 2]{B98}, we have uniform bounds on $\|\nabla w\|_{L^{\infty}(B_1^+)}$
\be \label{borsuk 1}
\|\nabla w\|_{L^{\infty}(B_1^+)}\le C(\mu, [A_{+}]_{C^{\alpha}(B_2^+)},N) \Big [ \|w\|_{L^{\infty}(B_2^+)} + \|\f\|_{C^{1,\alpha}(B_2')}  \Big ].
\ee
From maximum principle, we have $\|w\|_{L^{\infty}(B_2^+)} = \| w\|_{L^{\infty}(\partial B_2^+)}= \| \f \|_{L^{\infty}(\partial B_2^+)}\le M$ and from \ref{P1}, $\|\f\|_{C^{1,\alpha}(B_2')} \le M$. Plugging this information in \eqref{borsuk 1} we have 
$$
\|\nabla w\|_{L^{\infty}(B_1^+)}\le C(\mu, [A_{+}]_{C^{\alpha}(B_2^+)},N)M
$$
and then using the equation above in \eqref{borsuk 2}, 
\be\label{u+}
|u^+(x)| \le C(\mu, [A_{+}]_{C^{\alpha}(B_2^+)},N)M|x|, \qquad x\in B_1^+
\ee
Analogously, we have
\be\label{u-}
|u^-(x)| \le C(\mu, [A_{-}]_{C^{\alpha}(B_2^+)},N)M|x|,\qquad x\in B_1^+.
\ee
We add \eqref{u+} and \eqref{u-}, we obtain 
$$
u^+(x)+u^-(x) = |u(x)| \le C(\mu, [A_{\pm}]_{C^{\alpha}(B_2^+)},N)M|x|,\qquad \forall\; x\in B_1^+ .
$$
\end{proof}

\begin{remark}\label{linear growth for blowup}
We can check that for every $u\in \MM_1$ then $u_r\in \MM_r$ $($c.f. Remark \ref{rescaling invariance}$)$. Moreover, $u_r^{\pm}$ are $A^r_{\pm}$-subharmonic (c.f. Lemma \ref{subharmonic}, Lemma \ref{holder}) and $u_r$ satisfies \eqref{lg} in $B_{1/r}^+$. That is 
$$
|u_r(x)| \le C(\mu, [A_{\pm}]_{C^{\alpha}(B_2^+)},N)M\,|x|, \,\,x\in B_{1/r}^+.
$$
\end{remark}

\subsection{Energy estimates and compactness}
In this section, we prove the compactness lemma. That is for a given sequence $v_j \in \MM_{1}$ and $r_j\to 0^+$, the blow-ups of $v_j$ defined as $u_j:= (v_j)_{r_j} \in \MM_{r_j}$ (c.f. \eqref{blowup})  converge to $u_0\in \MM_{\infty}$ up-to a subsequence (in appropriate topologies c.f. Lemma \ref{blowup convergence}). For this purpose, we establish some uniform estimates for $\|u_j\|_{H^1(B_R^+)}$ which in turn allow us to use compact embeddings. 

Before proving energy estimates, we quote the following useful lemma 
\begin{lemma}\label{useful lemma} $($c.f. \cite[Lemma 6.1]{eg05}$)$
Let $Z(t)$ be a bounded non negative function defined in $[\rho , r]$. And assume that we have for $\rho \le s <t \le r$
$$
Z(s) \le \theta Z(t) + \frac{A}{|s-t|^2}+C
$$
for some $0\le \theta <1$ and $A,C\ge 0$. Then we have 
$$
Z(\rho) \le C(\theta) \Big [ \frac{A}{|\rho -r|^2 }+C \Big ].
$$
\end{lemma}

\begin{proposition}[Energy Estimates]\label{Energy Estimates}
Given $R>0$, $r>0$ such that $R\le \frac{1}{r}$, then for any \newline $u\in \MM_{1} $ we have 
\be\label{energy estimates}
\int_{B_R^+}|\nabla u_r|^2 \,dx\le C (N, \lambda_{\pm}, \mu, q, M,R).
\ee
\end{proposition}
{The proof we present below is completely different from \cite[Lemma 3.7]{MS22}, where we obtain bounds by controlling the $\L{2}{\nabla u_r}{B_R^+}$ by $L^2(B_R^+)$ norm of gradient of harmonic replacement of $u_r$ in $B_{2R}^+$. Instead, we present here a variational proof of Proposition \ref{Energy Estimates} which involve the use of Widman's hole filling technique. The presented proof also works in the context of \cite[Lemma 3.7]{MS22}, and other general scenarios in addition to the one in this paper.}
\begin{proof}
Suppose $s,t$ be such that $R\le t<s\le 2R$. Consider $\eta \in C_c^{\infty}(\R^N)$ be the following cutoff function  
$$
\eta := \begin{cases}
1\mbox{ in $B_t$}\\
0 \mbox{ in $\R^N\setminus B_s$}.
\end{cases}
$$
The cutoff function $\eta$ defined above can be taken in such a way that $0\le \eta\le 1$ and $|\nabla \eta | \le \frac{C(N)}{|s-t|}$. We consider the following test function $v\in H^1(B_s^+)$
$$
v:=u_r-\eta(u_r-\f_r) = u_r(1-\eta) + \eta \f_r.
$$
We can easily check that $v\in H^1(B_s^+)$. Moreover, $u_r,\eta$ and $\f_r$ are continuous functions (c.f. Lemma \ref{holder}), therefore $v-u_r \in C(B_s^+)$ and we can check that $u_r-v=0$ on $\partial B_s^+$. Since, $\partial B_s^+$ is Lipschitz, we have $v-u_r\in H_0^1(B_s^+)$ (c.f. \cite[Theorem 18.7]{GL09}). Hence $v$ is an admissible test function to compare minimality of $u_r$ for the functional $J(A^r_{\pm}, Q^r, \lambda_{\pm}, B_s^+)$. That is 
$$
\int_{B_s^+}\Big (  \< A^r(x,u_r) \nabla u_r,\nabla u_r  \>+ Q^r(x)\lambda (u_r)\Big )\,dx \le \int_{B_s^+} \Big ( \< A^r(x,v) \nabla v,\nabla v  \>+ Q^r(x)\lambda (v) \Big )\,dx.
$$
Now, we use ellipticity of $A^r_{\pm}$ (c.f. \ref{P2}) and boundedness of $Q^r(x)$ (c.f. \ref{P3}). Thus we have 
$$
\int_{B_s^+} |\nabla u_r|^2 \,dx  \le  C(\mu) \int_{B_s^+} |\nabla v|^2\,dx + C(\mu,\lambda_{\pm},q,N)s^N.
$$
Using the fact that $t<s$, we get
\be\label{part 0}
\begin{split}
\int_{B_t^+} |\nabla u_r|^2\,dx  &\le  C(\mu) \int_{B_s^+} |\nabla v|^2\,dx + C(\mu,\lambda_{\pm},q,N)s^N\\
&= C(\mu) \int_{B_s^+}  |\nabla (u_r(1-\eta) + \eta \f_r)|^2 \,dx +C(\mu,\lambda_{\pm},q,N)s^N\\
&\le C(\mu) \int_{B_s^+}\Big ( |\nabla (1-\eta)u_r|^2 + |\nabla \eta \f_r|^2\Big ) \,dx +C(\mu,\lambda_{\pm},q,N)s^N\\
&\le C(\mu)\int_{B_s^+} \Big ((1-\eta)^2 |\nabla u_r|^2 + |\nabla \eta|^2 |u_r|^2+ |\nabla \eta|^2 |\f_r|^2 + | \eta|^2 |\nabla \f_r|^2 \Big ) \,dx +C(\mu,\lambda_{\pm},q,N) s^N.
\end{split}
\ee
Since $\eta=1$ in $B_t^+$ and $0\le \eta \le 1$, we have $\int_{B_s^+} (1-\eta)^2 |\nabla u_r|^2 \,dx \le \int_{B_s^+\setminus B_t^+} |\nabla u_r|^2\,dx$. We also have $|\nabla \eta|\le \frac{C(N)}{|s-t|}$. From \ref{P1} we have $\Linfty{\f_r}{B_s^+}\le M$ and $\Linfty{\nabla \f_r}{B_s^+} \le M$. 
With all this information, we continue our calculations in \eqref{part 0}.
\be\label{part 1}
\begin{split}
\int_{B_t^+} |\nabla u_r|^2\,dx  & \le C_1(\mu) \int_{B_s^+\setminus B_t^+} |\nabla u_r|^2\,dx + \frac{C_2(\mu,M,N)}{|s-t|^2} \int_{B_s^+} (|u_r|^2 + |\f_r|^2)\,dx +C_3(\mu,\lambda_{\pm},q,N,M)s^N. 
\end{split}
\ee
Now, we add $C_1 \int_{B_t^+} |\nabla u_r|^2\,dx$ on both sides of \eqref{part 1} and we obtain the following
$$
\int_{B_t^+} |\nabla u_r|^2\,dx \le \frac{C_1}{C_1+1}\int_{B_s^+} |\nabla u_r|^2\,dx+ \frac{C_2'}{|s-t|^2} \int_{B_s^+} (|u_r|^2 + |\f_r|^2)\,dx +C_3's^N
$$
and since $s\le r$, we have 
\be\label{part 2}
\int_{B_t^+} |\nabla u_r|^2\,dx \le \frac{C_1}{C_1+1}\int_{B_s^+} |\nabla u_r|^2\,dx+ \frac{C_2'}{|s-t|^2} \int_{B_r^+} (|u_r|^2 + |\f_r|^2)\,dx +C_3' r^N.
\ee
From Lemma \ref{useful lemma} we obtain 
$$
\int_{B_R^+} |\nabla u_r|^2\,dx \le C(\mu,\lambda_{\pm},q,N,M) \Big [  \frac{1}{R^2} \int_{B_{2R}^+} (|u_r|^2 + |\f_r|^2)\,dx +R^N \Big ].
$$
Now, from Remark \ref{linear growth for blowup} we know that $|u_r(x)| \le C(\mu) |x| \le C(\mu)R$ in $B_R^+$ and from \ref{P1} we have $|\f_r| \le M|x|\le MR$ in $B_R^+$. Hence, after using this information in the above equation
\[
\begin{split}
\int_{B_R^+} |\nabla u_r|^2\,dx \le C(\mu,\lambda,q,N,M)  \Big [  \frac{1}{R^2} \int_{B_R^+} (|u_r|^2 + |\f_r|^2)\,dx +R^N \Big ] \le C(\mu,\lambda_{\pm},q,N,M)  R^N.
\end{split}
\]

\end{proof}


Now we prove an important result on convergence of blowup of minimizers, which is crucial in proving Theorem \ref{tangential touch}.

\begin{proposition}[Compactness]\label{blowup convergence}
Let $v_j \in \MM_1$ be a sequence and $r_j\to 0^+$. We denote $u_j := (v_j)_{r_j}$. Then, for every $R>0$ $($such that $R<\frac{1}{r_j}$ for all $j\in \N$$)$ the sequence $u_j$ is pre-compact in weak-$H^1(B_R^+)$ and $L^{\infty}(B_R^+)$ for every $R>0$. 
Moreover, if $u_0$ is a subsequential limit of $u_j$ in above mentioned topologies, then $u_0\in \MM_{\infty}$.
\end{proposition}

\begin{proof}
By change of variables, we know that $u_j\in \MM_{r_j}$. We set the notation for the functional $J_j(\cdot;B_R^+)$ 
\be\label{Jj}
J_j(v;B_R^+):=\int_{B_{R}^+} \Big ( \< A_j(x,v) \nabla v,\nabla v  \>+ Q_j(x)\Lambda (v) \Big )\,dx\qquad  v\in H_{\phi_j}^1(B_{2/r_j}^+).
\ee
Here $ A_{j,\pm}(x)$, $\phi_j(x)$, $\lambda_{\pm}$ and $Q_j(x)$  satisfy \ref{P1}-\ref{P5} for $r=r_j$.

We can see that $J_j$ satisfy the structural conditions \cite[Theorem 7.3]{eg05}. From Lemma \ref{holder}, we have $\|u_j\|_{C^{\alpha_0}(\overline{B_R^+})} \le C(M,\lambda_{\pm},\mu)$. Hence the sequence $u_j$ is uniformly bounded and equicontinuous in $\overline{B_R^+}$. By Arzela-Ascoli theorem, there exists $u_0 \in C^{0,\alpha_0}(B_R^+)$ such that 
\be\label{uniform convergence}
 u_j \to u_0 \qquad\mbox{ in $L^{\infty}({B_R^+})$ up-to a subsequence}.
\ee
Also, from Lemma \ref{Energy Estimates} for all $j\in \N$, we have 
$$
\int_{B_R^+}|\nabla u_j|^2\,dx \le C (N, \lambda_{\pm}, \mu, q, M,R).
$$
Moreover, the linear growth condition  is preserved under blowup (c.f. Remark \ref{linear growth for blowup}), that is 
\be\label{lg for seq}
|u_{j}(x)| \le C(\mu)M |x| \mbox{ for all $x\in B_R^+$}
\ee
then we have 
$$
\int_{B_R^+}|u_j|^2\,dx \le {C(\mu)M} |B_R^+|.
$$
Therefore, $\|u_j\|_{H^1(B_R^+)}$ is a bounded sequence and $\{u_j\}$ converges in weak-$H^1(B_R^+)$ to $u_0$, up to a subsequence (which we again rename as $u_j$). By assumptions \ref{P1}-\ref{P5} on $A_j$, $Q_j$ and $\phi_j$ (for $r=r_j$), we can see that there exists constants $a_{\pm} \in (\mu, \frac{1}{\mu})$ and $q\le q_0 \le \frac{1}{q}$ such that in the limit ${j\to \infty}$
\be\label{convergence}
\begin{split}
A_{\pm,j}(x) \to a_{\pm}Id \qquad &\mbox{ uniformly in $B_R^+$}\\
Q_j(x) \to q_0 \qquad &\mbox{ uniformly in $B_R^+$}\\
\f_j \to 0 \qquad &\mbox{ uniformly in $B_R'$}.
\end{split}
\ee
Moreover, since $u_0\in C(\overline {B_R^+})$ and $u_0 = \lim_{j\to \infty}u_j= \lim_{j\to \infty} \f_j=0$ on $B_R'$, $u_0$ satisfies \ref{G2} and \ref{G3} on $B_R^+$ for all $R>0$. By taking $\lim_{j\to \infty}$ in \eqref{lg for seq} we obtain $|u_0(x)| \le C(\mu, M)|x|$ for all $x\in B_R^+$. This proves \ref{G1} in $B_R^+$ for all $R>0$.

In order to verify \ref{G4}, we claim that $u_0$ is minimizer of $J_0(\cdot;B_R^+)$ in the set $H_{u_0}^1(B_R^+)$. Here $J_0(\cdot;B_R^+)$ is denoted by the functional 
\be\label{J0}
J_0(v;B_R^+) := \int_{B_R^+}\Big ( a(v)|\nabla v|^2+ q_0\Lambda(v) \Big )\,dx, \qquad v\in H^1(B_R^+).
\ee
First let us show that 
\be \label{liminf}
J_0(u_0;B_R^+) \le \liminf_{j\to \infty} J_j(u_j;B_R^+)
\ee
We look at the right hand side of \eqref{liminf} term by term. We claim that
\be\label{last term}
\begin{split}
\int_{B_R^+}q_0 \Lambda(u_0) \,dx &\le \liminf_{j\to \infty} \int_{B_R^+} Q_j(x) \Lambda(u_j)\,dx.\\
\end{split}
\ee 
To prove the claim \ref{last term}, we first show that for almost every $x\in B_R^+$, we have
\be\label{char}
q_0\lambda_{+}\chi_{\{u_0>0\}}(x)+ q_0\lambda_- \chi_{\{u_0\le0\}}(x) \le \liminf_{j\to \infty} \left ( Q_j(x)\lambda_+ \chi_{\{u_j>0\}}(x)+Q_j(x)\lambda_-\chi_{\{u_j\le0\}}(x)\right ).
\ee
Let $x_0\in B_R^+ \cap \big ( \{u_0>0\} \cup \{u_0<0\} \big )$, then from the uniform convergence of $u_j$ to $ u_0$, we can easily see that $u_j(x_0)$ attains the sign of $u_0(x_0)$ for sufficiently large value of $j$. Hence, combining this information with \eqref{convergence}, the claim\eqref{char} holds true in the set $\big ( \{u_0>0 \} \cup \{u_0<0\} \big ) \cap B_R^+$.

If $x_0\in \{u_0=0\}$. Then left hand side of \eqref{char} is equal to 
$$
q_0\lambda_{+}\chi_{\{u_0>0\}}(x_0)+ q_0\lambda_- \chi_{\{u_0\le0\}}(x_0) = q_0\lambda_{-}.
$$
We observe that in the set $\{u_0 = 0\}\cap B_R^+$,  the right hand side of \eqref{char} is the following
$$
Q_j(x_0)\lambda_+ \chi_{\{u_j>0\}}(x_0)+Q_j(x_0)\lambda_-\chi_{\{u_j\le0\}}(x_0) = 
\begin{cases}
Q_j(x_0) \lambda_+, \qquad \mbox{if $u_j(x_0)>0$}\\
Q_j(x_0)  \lambda_-, \qquad \mbox{if $u_j(x_0)\le 0$}.
\end{cases}
$$
Since $\lambda_- <\lambda_+$ (c.f \ref{P3}), the right hand side in \eqref{char} is always greater than or equal to $ Q_j(x_0) \lambda_- $. That is we have  for $x_0\in \{u_0=0\}\cap B_R^+$ 
\[
\begin{split}
q_0\lambda_{+}\chi_{\{u_0>0\}}(x_0)+ q_0\lambda_- \chi_{\{u_0\le0\}}(x_0) &= q_0\lambda_{-}\\
&= \liminf_{j\to \infty} Q_j(x_0) \lambda_-\\
&\le \liminf_{j\to \infty} \left ( Q_j(x_0)\lambda_+ \chi_{\{u_j>0\}}(x_0)+Q_j(x_0)\lambda_-\chi_{\{u_j\le0\}}(x_0) \right ).
\end{split}
\]
Thus, \eqref{char} is proven for all $x\in B_R^+$ and hence \eqref{last term} holds by Fatou's lemma. 

Now, recalling that $A_{\pm}(0)=a_{\pm} Id$ we see that 
\be\label{first term}
\int_{B_R^+}\< A_{j,\pm}(x) \nabla u_j^{\pm},\nabla u_j ^{\pm}\>\,dx= \int_{B_R^+}\<  (A_{j,\pm}(x)-a_{\pm} Id) \nabla u_j^{\pm},\nabla u_j ^{\pm}\>\,dx + \int_{B_R^+} a_{\pm} |\nabla u_j^{\pm}|^2\,dx
\ee
from Lemma  \ref{Energy Estimates}, we know $\int_{B_R^+}|\nabla u_j^{\pm}|^2\,dx$ is uniformly bounded for all $j\in \N$ and $A_{j,\pm} \to a_{\pm}$ uniformly as $j\to \infty$ (c.f. \eqref{convergence}), therefore, the first term on the right hand side of \eqref{first term} converges to zero as $j\to \infty$. Regarding the second part, we observe that $u_j\rightharpoonup u_0$ weakly in $H^1(B_R^+)$.  Hence, from {\cite[Proposition 3.7]{MS21}} $u_j^{\pm} \rightharpoonup u_0^{\pm}$ weakly in $H^1(B_R^+)$.  By the lower semi-continuity of $H^1$ norm in the weak topology and \eqref{first term}, we obtain,
\be\label{wlsc}
\int_{B_R^+} a_{\pm}^2 |\nabla u_0^{\pm}|^2 \le  \liminf_{j\to \infty} \int_{B_R^+} a_{\pm}^2 |\nabla u_j^{\pm}|^2\,dx= \liminf_{j\to \infty} \int_{B_R^+}\< A_{\pm}(x) \nabla u_j^{\pm},\nabla u_j ^{\pm}\>\,dx.
\ee
By adding \eqref{wlsc}, \eqref{last term} and applying \cite[Theorem 3.127]{D17} to the summation, we obtain \eqref{liminf}. In order to prove that the function $u_0\in \MM_{\infty}$, it only remains to verify that $u_0$ satisfies \ref{G4} in $B_R^+$ for all $R>0$. That is $u_0$ is a minimizer of $J_0(\cdot; B_R^+)$ for all $R>0$. In this direction, we consider $w\in H_{u_0}^1(B_R^+)$. We make use of \eqref{liminf} to prove our claim. Let us construct a competitor for the minimality of $u_j$ for functional $J_j$.  We define two cutoff functions $\eta_{\delta}: \R^N\to \R$ and $\theta :\R \to \R$ as follows,
$$
\eta_{\delta}(x):=
\begin{cases}
1,\; x\in B_{R-\delta}\\
0, \;x \in \R^N\setminus B_{R}
\end{cases},
\theta(t):=
\begin{cases}
1,\;|t|\le 1/2\\
0,\; |t| \ge 1.
\end{cases}
$$
we can take $|\nabla \eta_{\delta}|\le \frac{C(N)}{\delta}$.  We define $\theta_j(x)= \theta(\frac{x_N}{d_j})$, for a sequence $d_j\to 0$, which we choose in later steps of the proof. Let $w_j^{\delta}$ be a test function defined as 
\be\label{test function}
w_j^{\delta} := w + (1-\eta_{\delta}) (u_j-u_0) +\eta_{\delta} \theta_{j} \f_j.
\ee
Since, the function $w_j^{\delta} -w = (1-\eta_{\delta}) (u_j-u_0)+\eta_{\delta} \theta_{j} \f_j$ is continuous in $\overline{B_R^+}$ and is pointwise equal to zero on $\partial B_R^+$ (because $u_0=0$ on $B_R'$, c.f. \eqref{convergence}). Since $\partial B_R^+$ is a Lipschitz surface in $\R^N$,  $u_j-w_j^{\delta} \in H_0^1(B_R^+)$. See that from minimality of $u_j$, we have 
$$
\int_{B_R^+} \Big (\<  A_j(x, u_j) \nabla u_j, \nabla u_j \> + Q_j(x)\Lambda(u_j) \Big ) \,dx \le  \int_{B_R^+} \Big ( \<  A_j(x, w_j^{\delta}) \nabla w_j^{\delta}, \nabla w_j^{\delta} \> +Q_j(x) \Lambda(w_j^{\delta}) \Big ) \,dx.
$$
From \eqref{liminf} and \cite[Theorem 3.127]{D17}, we obtain
\be\label{minimality}
\begin{split}
J_0(u_0;B_R^+)  &\le \liminf_{j\to \infty}  \int_{B_R^+}\Big ( \<  A_j(x, w_j^{\delta}) \nabla w_j^{\delta}, \nabla w_j^{\delta} \> +Q_j(x) \Lambda(w_j^{\delta}) \Big )\,dx\\
&\le \limsup_{j\to \infty}  \int_{B_R^+}\Big ( \<  A_j(x, w_j^{\delta}) \nabla w_j^{\delta}, \nabla w_j^{\delta} \> +Q_j(x) \Lambda(w_j^{\delta}) \Big )\,dx\\
&\le \limsup_{j\to \infty}  \Big (  \int_{B_R^+}\Big ( \<  A_j(x, w_j^{\delta}) \nabla w_j^{\delta}, \nabla w_j^{\delta} \> \,dx \Big ) + \limsup_{j\to \infty} \Big (\int_{B_R^+}  Q_j(x) \Lambda(w_j^{\delta}) \,dx\Big )
\end{split}
\ee
Since, from \eqref{convergence} that $A_{j,\pm} \to a_{\pm}$ uniformly as $j\to \infty$, following the same reasoning as in \eqref{first term}, we have
\[
\begin{split}
\limsup_{j\to \infty} \int_{B_R^+} \<  A_j(x, w_j^{\delta}) \nabla w_j^{\delta}, \nabla w_j^{\delta} \>\,dx &= \limsup_{j\to \infty}  \int_{B_R^+} \<  a(w_j^{\delta}) \nabla w_j^{\delta}, \nabla w_j^{\delta} \>\,dx 
\end{split}
\]
Our claim is that after passing the limit $\delta\to 0$, the right hand side of \eqref{minimality} is in-fact $J_0(w;B_R^+)$. To prove this, we argue exactly the way as in the proof of \cite[Lemma 3.8]{MS22}. In \cite[Lemma 3.8]{MS22}, it is shown that by choosing $d_j\to 0$ such that $\frac{r_j^{\alpha}}{d_j}\to 0$ we can prove that as $\delta\to 0$ and $j\to \infty$, $\int_{B_R^+} |\nabla (w_j^{\delta}-w)|^2\,dx\to 0$, in other words 
\be\label{eq14}
\lim _{\delta \to 0}\left ( \limsup_{j\to \infty } \int_{B_R^+} | \nabla w_j^{\delta}|^2 \,dx \right ) = \int_{B_R^+} |\nabla w|^2\,dx
\ee
and therefore from \cite[Proposition 3.7]{MS21}
\be\label{eq15}
\lim _{\delta \to 0}\left ( \limsup_{j\to \infty } \int_{B_R^+} a_{\pm} |\nabla (w_j^{\delta})^{\pm}|^2\,dx \right ) = \int_{B_R^+}a_{\pm} |\nabla w^{\pm}|^2\,dx.
\ee
Following the steps as in the proof of \cite[Lemma 3.8]{MS22}, it can also be shown that we have 
\be\label{eq16}
\lim_{\delta \to 0} \left ( \limsup_{j\to \infty} \int_{B_R^+}Q_j(x) \Lambda(w_j^{\delta}) \,dx \right )=  \int_{B_R^+}q_0 \Lambda(w) \,dx.
\ee
We now use \eqref{eq15} and \eqref{eq16} in \eqref{minimality}, and obtain $J_0(u_0;B_R^+) \le J_0(w;B_R^+)$ for all $w\in H_{u_0}^1(B_R^+)$ for all $R>0$. This verifies \ref{G4} for $u_0$. Hence we conclude the proof of Proposition \ref{blowup convergence}.
\end{proof}

\subsection{Non-degeneracy results}
In this section we prove various results on non-degenerate behavour of $u\in \MM_1$ near their free boundaries and also near contact point. The following proposition adapts ideas from \cite[Theorem 3.1]{ACF84}. Although, the proof of the following lemma is exactly the same as in \cite[Proposition 3.9]{MS22}, we make small changes to fit in the context of transmission problems.

\begin{proposition}[Non-degeneracy near the free boundary] \label{non degeneracy}
For $u\in \MM_{r_0}$,  and $x_0\in B_{2/r_0}^+$. For every $0<\kappa<1$ there exists a constant $c>0$, which depends only on $\k,\lambda_{\pm},\mu,q$ such that for all $B_r(x_0)\subset B_2^+$, we have 
\be
\frac{1}{r} \fint _{\partial B_r(x_0)} u^+\,d\HH^{N-1}(x)\,<\,c(\k,\lambda_{\pm},\mu,q,N) \implies \mbox{$u^+=0$ in $B_{\kappa r}(x_0)$}.
\ee 
\end{proposition}

\begin{proof}
Let $\gamma = \frac{1}{r} \fint _{B_r(x_0)} u^+\,dx$. Since, by elliptic regularity theory, $u$ is locally $C^{1,\alpha}$ in $\{u>0\}$. Then, for almost every $\epsilon\in (0, \Linfty{u^+}{B_r(x_0)})$, $B_r\cap \partial \{u>\eps\}$ is a $C^{1,\alpha}$ surface. For one such small $\eps>0$, we consider the test function $v_{\eps}$ given by
$$
\begin{cases}
\dive(A_+(x)\nabla v_{\eps}) =0 \qquad \;\;\; \mbox{in $(B_r(x_0)\setminus B_{\kappa r}(x_0))\cap \{u>\eps\}$}\\
v_{\eps} = u\qquad \qquad \;\; \qquad \qquad \mbox{in $B_{r}(x_0)\cap \{u\le \eps \}$}\\
v_{\eps}=\eps \qquad \;\; \qquad \qquad \qquad\mbox{in $B_{\k r}(x_0)\cap \{u>\eps\}$}\\
v_{\eps}=u \qquad \; \; \qquad \qquad \qquad \mbox{on $\partial B_r(x_0)$}
\end{cases}
$$ 
we know $v_{\eps}\in H^1(B_r(x_0))$, thanks to \cite[Theorem 3.44]{DD12} and \cite[Theorem 4.6]{GE15}. To ensure the existence of limit $\lim_{\eps\to 0}v_{\eps}$ exists in weak sense in $H^1(B_r(x_0))$ and strong sense in $L^2(B_r(x_0))$ (up-to a subsequence), we show that $ v_{\eps}$ is bounded in $H^{1}(B_r(x_0))$. To prove this, let $G$ be the Green function for $L(v)=\dive(A_+(x)\nabla v)$ in the ring $B_r(x_0)\setminus B_{\k r}(x_0)$. Then, for a function $w$ such that 
$$
\begin{cases}
\dive(A_+(x)\nabla w)=0 \; &\mbox{in $B_r(x_0)\setminus B_{\k r}(x_0)$}\\
w=u \; &\mbox{on $\partial B_r(x_0)\cap \{u>\eps\}$}\\
w=\eps \; &\mbox{elsewhere on $\partial (B_r(x_0)\setminus B_{\k r}(x_0))$}.
\end{cases}
$$
We apply \cite[Theorem 3.3 (vi)]{GW82}, for $\b{x}\in \partial B_{\k r}(x_0)$ and for a sequence $\{x_k\} \subset B_r(x_0)\setminus B_{\k r}(x_0)$ such that $x_k \to \bar x \in \partial B_{\k r}(x_0)$.
\be\label{green estimates}
\begin{split}
\big | \nabla w(\bar x)\big | &\le {C(N,\mu)} \lim_{k\to \infty} \int_{\partial B_r(x_0)\cap \{u>\eps\}}\Big | \nabla_x \left ( \frac{\partial}{\partial _{\nu_y}}G(x_k,y) \right ) (u-\eps)^+ \Big |\,dx\\
& \le {C(N,\mu)} \lim_{k\to \infty} \int_{\partial B_r(x_0)} \frac{1}{|x_k-y|^N}(u-\eps)^+\,d\HH^{N-1}(y)\\
&\le \frac{C(N,\mu, \alpha)}{(1-\k)^N}\frac{1}{r}\fint _{\partial B_r} (u-\eps)^+\,d\HH^{N-1}(y)\le C(\mu,N,\k) \gamma\;\mbox{on $\partial B_{\k r}(x_0)$}.
\end{split}
\ee
We can easily check that $w\ge v_{\eps}$ on $B_r(x_0)\setminus B_{\k r}(x_0)$.  By comparison principle, $w\ge v_{\eps}$ in $B_r(x_0)\setminus B_{\k r}(x_0)$. Since $w=v_{\eps}$ on $\partial B_{\k r}(x_0)$, we have
\be\label{eq 2.6}
|\nabla v_{\eps}|\le |\nabla w|\le C(\k) \gamma\;\;\mbox{on $\partial B_{\k r}(x_0)$}.
\ee
Let us denote $D_{\eps}:= (B_r(x_0)\setminus B_{\k r}(x_0))\cap \{u>\eps\}$. Since $\dive (A_+(x)\nabla v_{\eps})=0$ in $D_{\eps}$, we have by divergence theorem
\[
\begin{split}
\int_{D_{\eps}} (A_+(x)\nabla v_{\eps})\cdot \nabla (v_\eps-u)\,dx&= \int_{\partial B_{\k r}(x_0)\cap \{u>\eps\}} (u-v_{\eps})(A_+(x)\nabla v_{\eps})\cdot \nu(y)\, d\HH^{N-1}(y)\\
&\le \frac{1}{\mu} \int_{\partial B_{\k r}(x_0) \cap \{ u>\eps \}} |u-\eps||\nabla v_{\eps}|\,d\HH^{N-1}(y)\\
&\le \frac{C (N,\k) \gamma}{\mu}\int_{\partial B_{\k r}(x_0) \cap \{ u>\eps \}} |u-\eps|\,d\HH^{N-1}(y)\le C_0(u).
\end{split}
\]
One can refer to \cite[equation (3.4)]{ACF84} for a justification of use of divergence theorem in $D_{\eps}$ which is just a Lipschitz domain. In the last line of calculations above, we have used \eqref{eq 2.6}. From the above calculations, we can write 

\begin{align*}
\;\; \;\;& \int_{D_{\eps}} (A_+(x)\nabla v_{\eps})\cdot \nabla (v_\eps-u) \>\,dx \le C_0(u)\\
 \Rightarrow&  \int_{D_{\eps}} (A_+(x)\nabla v_{\eps})\cdot \nabla v_\eps\,dx  \le C_0(u)+\int_{D_{\eps}} (A_+(x)\nabla v_{\eps})\cdot \nabla u\,dx\\
\Rightarrow &\mu \int_{D_{\eps}}|\nabla v_{\eps}|^2\,dx \le C_0(u) + \frac{1}{\mu} \int_{D_{\eps}} |\nabla v_{\eps}||\nabla u|\,dx\\
\Rightarrow & \mu \int_{D_{\eps}}|\nabla v_{\eps}|^2\,dx\le C_0(u)+\frac{\eps_0}{2\mu} \int_{D_{\eps}} |\nabla v_{\eps}|^2\,dx+\frac{1}{2\eps_0 \mu}\int_{D_\eps} |\nabla u|^2\,dx
\end{align*}

putting very small $\eps_0>0$ in the last inequality, we have 
$$
\int_{D_{\eps}}|\nabla v_{\eps}|^2\,dx \le C_1(u)
$$
and thus, since $v_{\eps}=\eps \Rightarrow \nabla v_{\eps}=0$ in $B_{\k r}(x_0)\cap \{u>\eps\}$ and $v_\eps = u $ in $B_{r}(x_0) \setminus D_{\eps}$, 
$$
\int_{B_r(x_0)}|\nabla v_\eps|^2\,dx = \int_{D_{\eps}}|\nabla v_{\eps}|^2\,dx+\int_{B_r(x_0)\setminus D_{\eps}}|\nabla u|^2\,dx\le\,C_2(u).
$$
We note that by definition of $v_{\eps}$ and comparison principle, we have $u^-\le v_{\eps} \le u$.  Therefore 
$$
\int_{B_r(x_0)}| v_\eps|^2\,dx \le \int_{B_r(x_0)}| u|^2\,dx.
$$
This proves $v_{\eps}$ is uniformly bounded in $H^1(B_r(x_0))$.  Therefore there exists a limit $v = \lim_{\eps \to 0} v_{\eps} $ in weak $H^1$ sense, such that $v$ satisfies the following
\be\label{prop of v}
\begin{cases}
\dive(A_+(x)\nabla v) =0 \qquad \mbox{in $(B_r(x_0)\setminus B_{\k r}(x_0))\cap \{u>0\}$}\\
v = u\qquad \qquad \;\; \qquad \qquad \mbox{in $B_{r}(x_0)\cap \{u\le 0 \}$}\\
v =0 \qquad \;\; \qquad \qquad \qquad\mbox{in $B_{\k r}(x_0)\cap \{u>0\}$}\\
v=u \qquad \; \; \qquad \qquad \qquad \mbox{on $\partial B_r(x_0)$}.
\end{cases}
\ee
The above properties of the function $v\in H^1(B_r(x_0))$ are verified in \cite[Proposition 3.9]{MS22}. Let us use the function $v$ as a test function with respect to minimality condition on $u$ in $B_r(x_0)$, we have
\[
\begin{split}
\int_{B_r(x_0)}\Big ( \< A(x,u)\nabla u,\nabla u\>+ Q(x) \Lambda (u) \Big ) \,dx& \le \int_{B_r(x_0)}\Big ( \< A(x,v)\nabla v,\nabla v\> + Q(x)\Lambda (v) \Big )\,dx\\
\end{split}
\]
we observe that, $v=u$ in $\{u\le 0\}$ and $\{ v>0 \} \subset \{u>0\}$. Thus, the integration in the set $\{u\le 0\}\cap B_r(x_0)$ gets cancelled from both sides and we are left with the equations below. Set $D_0:=B_r(x_0)\setminus B_{\k r}(x_0)\cap \{u>0\}$, since $q\le Q \le \frac{1}{q}$ we have
\[
\begin{split}
\int_{B_r(x_0)\cap \{u> 0\}} \Big (  \< A(x) \nabla u,\nabla u \> - \< A(x) \nabla v,\nabla v \> \Big )\,dx &\le \int_{B_{ r}(x_0)\cap \{u>0\}} Q(x)(\Lambda(v)-\Lambda(u))\,dx\\
&=  \frac{\lambda_0}{q} |B_{\k r}(x_0)\cap \{u>0\}| . \qquad(\lambda_0: = -(\lambda_ + - \lambda_-)). \\ 
\end{split}
\] have second equality above because $\chi_{\{v>0\}}= \chi_{u>0}$ in $D_0$. Since $v=0$ in $D_0$, we have 
$$
\int_{B_{\k r}(x_0)\cap \{u>0\}} \< A_+(x) \nabla u,\nabla u \>\,dx + \int_{D_0} \Big ( \< A_+(x) \nabla u,\nabla u \> - \< A_+(x) \nabla v,\nabla v \> \Big )\,dx\le \int_{B_{\k r(x_0)}\cap \{u>0\}} \frac{\lambda_0}{q} \,dx.
$$
Using the ellipticity of $A_+$ and shuffling the terms in the above equation, we obtain
\be\label{eq2.5}
\begin{split}
\int_{B_{\k r}(x_0)\cap \{u>0\}} \Big (\mu |\nabla u|^2-\frac{\lambda_0}{q} \Big ) \,dx &\le \int_{D_0}  \Big (\< A_+(x) \nabla v,\nabla v \>- \< A_+(x) \nabla u,\nabla u \> \Big ) \,dx\\
&= \int_{D_0} \< A_+(x) \nabla (v-u), \nabla (v+u) \>\,dx\\
&=  \int_{D_0} \< A_+(x) \nabla (v-u), \nabla (u-v+2v) \>\,dx\\
& \le 2\int_{D_0} \< A_+(x)\nabla v, \nabla (v-u) \>\,dx \\
& \le  \liminf_{\eps\to 0} 2\int_{D_0} \< A_+(x) \nabla v_{\eps}, \nabla (v_{\eps}-u) \>\,dx\\
&= \liminf_{\eps\to 0} 2 \int_{D_{\eps}} \< A_+(x) \nabla v_{\eps}, \nabla (v_{\eps}-u) \>\,dx\\
&=\liminf_{\eps\to 0} 2\int_{\partial B_{\k r}(x_0)\cap \{u>\eps\}} (u-\eps) ( A_+(x)\nabla v_{\eps})\cdot \nu \,dx\\
&\le \liminf_{\eps\to 0} \frac{2}{\tau} \int_{\partial B_{\k r}(x_0)\cap \{u>\eps\}}(u-\eps)\big | \nu\cdot \nabla v_{\eps} \big |\,dx\,:=\,M_0.
\end{split}
\ee
The second to last equality in above calculation is obtained from integration by parts, its justification can be found in \cite[equation (3.4)]{ACF84}. From \eqref{eq2.5} and \eqref{eq 2.6},  and using the trace inequality in $H^1(B_{\k r})$ we have (for some different constant $C(\k)$),

\be\label{eq2.8}
\begin{split}
M_0 &\le C(\mu,N,\k) \gamma \int_{\partial B_{\k r(x_0)}} u^+\,d\HH^{N-1}(x)\\
&\le C(\mu,N, \k) \gamma   \int_{B_{\k r(x_0)}} \Big (  |\nabla u^+|+\frac{1}{r} u^+ \Big )\,dx \\
& \le C(\mu,N, \k) \gamma \Bigg [ |B_{\k r(x_0)}\cap \{u>0\}|^{1/2}\left ( \int_{B_{\k r(x_0)}}|\nabla u^+|^2\,dx \right )^{1/2}+ \frac{1}{r}\sup_{B_{\k r}(x_0)} (u^+)\big | \{B_{\k r}(x_0)\cap \{u>0\}\}\big | \Bigg ]\\
& \le C(\mu,N, \k) \gamma \Bigg [  \frac{1}{2 \sqrt{-\lambda_0/q} } \int_{B_{\k r}(x_0)\cap \{u>0\}} |\nabla u^+|^2\,dx +  2 \sqrt{-\lambda_0/q}  |B_{\k r(x_0)}\cap \{u>0\}| \\ 
&\qquad \qquad  \qquad \qquad  \qquad \qquad \qquad \qquad  \qquad \qquad   \qquad \qquad \qquad + \frac{1}{r} \sup_{B_{\k r}(x_0)}(u^+) \int_{B_{\k r}(x_0)\cap \{u>0\}}1\,dx \Bigg ] \\
&= \frac{C(\mu,N, \k)\gamma }{2 \sqrt{-\lambda_0/q} }\left ( \int_{B_{\k r}(x_0)\cap \{u>0\}} |\nabla u^+|^2 - \frac{\lambda_0}{q}  \,dx \right )+\frac{C(\mu,N,\k)\gamma}{\frac{(-\lambda_0)}{q}  r}\sup_{B_{\k r}(x_0)}(u^+) \int_{B_{\k r}(x_0)\cap \{u>0\}}\frac{(-\lambda_0)}{q}   \,dx
\end{split}
\ee
we have used Hölder's inequality and then Young's inequality above. From \cite[Lemma 3.1]{MS21}, we know that $u^+$ is $A_+$-subharmonic in $B_r(x_0)$. If $G'$ is the Green's function for $L'(v)=\dive(A_+(x)\nabla v)$ in $B_r(x_0)$, then by comparison principle
$$
u^+(x)\le \int_{\partial B_r(x_0)} u^+(y)  \left  ( A_+(y)\nabla_y G'(x,y)\right )\cdot \nu_y \,d\HH^{N-1}(y)\;\; \forall x\in B_{\k r}(x_0).
$$
Since for all $y\in \partial B_r(x_0)$ and $x\in B_{\k r}(x_0)$, we have $\frac{1}{|x-y|^{N-1}}\le \frac{C(\k)}{r^{N-1}}$, then using the Green's function estimates c.f. \cite[Theorem 3.3 (v)]{GW82} we get 
\be\label{eq2.9}
\begin{split}
\sup_{B_{\k r}(x_0)}u^+&\le C(N,\mu) \int_{\partial B_r(x_0)} \frac{u^+(y)}{|x-y|^{N-1}}\,d\HH^{N-1}(y)\\
&\le C(N,\k,\mu) \fint_{\partial B_r} u^+\,d\HH^{N-1}(y)= C(N,\k,\mu)\gamma r.
\end{split}
\ee
Let us denote the integral $\int_{B_{\k r}(x_0)\cap \{u>0\}} \Big ( |\nabla u|^2-\frac{\lambda_0}{q} \Big ) \,dx$ by the letter $\II$. 
$$
\II:=\int_{B_{\k r}(x_0)\cap \{u>0\}} \Big ( |\nabla u|^2-\frac{\lambda_0}{q} \Big ) \,dx.
$$
By using, \eqref{eq2.5} and \eqref{eq2.9} in \eqref{eq2.8} and we have
\[
\begin{split}
\mu \II & \le 
\frac{C(\mu,\k,N)\gamma }{2 \sqrt{-\lambda_0/q} }\II +\frac{C(\mu,\k,N)\gamma}{(-\lambda_0/q)  r}\sup_{B_{\k r}(x_0)}(u^+) \int_{B_{\k r}(x_0)\cap \{u>0\}}\frac{(-\lambda_0)}{q} \,dx\\
&\le \frac{C(\mu,\k,N) \gamma}{\mu \sqrt{-\lambda_0} } \left (  1+ \frac{C(\k)\gamma}{\sqrt{-\lambda_0}} \right )\II.
\end{split}
\]
If $\gamma$ is sufficiently small, then $\II=0$  
$$
\frac{(-\lambda_0)}{q}|B_{\k r}(x_0)\cap \{u>0\}| \le \int_{B_{\k r}(x_0)\cap \{u>0\}}  \Big (|\nabla u|^2+\frac{(-\lambda_0)}{q} \Big ) \,dx=0
$$
in particular $|\{u>0\}\cap B_{\k r(x_0)}|=0$, that is $u^+=0$ almost everywhere in $B_{\k r}(x_0)$.
\end{proof}


\section{Proof of Theorem \ref{tangential touch}.}\label{main proof}

Before moving into the proof of Theorem \ref{tangential touch}, we need to show that the the positivity sets $\{u_j>0\} \cap B_R^+$ for the blowups $u_j$ as in Proposition \ref{blowup convergence}, converge in $L^1(B_R^+)$ to the positivity set of blowup limit $\{u_0>0\}\cap B_R^+$. 

\begin{lemma}
Let $u_0$ and $u_k$ be as in Proposition \ref{blowup convergence}. Then, for a subsequence of $u_k$, for any $R>0$ we have 
\be\label{pw convergence}
\chi_{\{u_k>0\} \cap B_R^+} \to \chi_{\{u_0>0\} \cap B_R^+}\qquad \mbox{ a.e. in $B_R^+$}.
\ee
This in turn implies 
\be\label{L1 convergence}
\chi_{\{u_k>0\} \cap B_R^+} \to \chi_{\{u_0>0\} \cap B_R^+}\qquad \mbox{ in $L^1(B_R^+)$}.
\ee
\end{lemma}

\begin{proof}
From Proposition \ref{blowup convergence}, we can consider a subsequence of $u_k$ such that $u_k\to u_0$ in $L^{\infty}(B_R^+)$.
Let $x\in B_R^+$. If $x\in \{u_0>0\}\cap B_R^+$ (or $\chi_{\{u_0>0\}\cap B_R^+}(x)=1$), then for  sufficiently large $k$, $u_k(x)$ attains the sign of $u_0(x)$. Thus we conclude that 
$$
\mbox{ $\chi_{\{u_k>0\}\cap B_R^+} (x) \to \chi_{\{u_0>0\}\cap B_R^+} (x)$ as $k\to \infty$ for all $x\in \{u_0>0\}\cap B_R^+$}. 
$$
If $x\in \{u_0\le 0\}^o \cap B_R^+$ (or $\chi_{\{u_0>0\}\cap B_R^+}(x)=0$), then there exists $\delta>0$ such that \newline $B_{\delta}(x) \subset \{u_0\le 0\} \cap B_R^+$. In other words, $\frac{1}{\delta} \fint_{\partial B_\delta(x)} u_0^+ \,d\HH^{N-1} =0$, by the uniform covergence of $u_k$ to $u_0$ in $B_R^+$ (c.f. Proposition \ref{blowup convergence}) we obtain 
\be\label{c+}
\frac{1}{\delta} \fint_{\partial B_\delta(x)} u_k^+ \,d\HH^{N-1} \le \frac{1}{2}c(\lambda_{\pm},\mu,q,N)\qquad\mbox{for $k$ sufficiently large.}
\ee
Here $c(\lambda_{\pm},\mu,q,N)$ is as in Proposition \ref{non degeneracy}. This implies $u_k \le 0$ in $B_{\frac{\delta}{2}}(x)$ (c.f. Proposition \ref{non degeneracy}). In particular, $\chi_{\{u_k(x)> 0\}}(x) =0$ for $k$ sufficiently large. This way,  
\be\label{c-}
\mbox{ $\chi_{\{u_k>0\}\cap B_R^+} (x) \to \chi_{\{u_0> 0\}\cap B_R^+} (x)$ as $k\to \infty$ for all $x\in \{u_0\le 0\}\cap B_R^+$}. 
\ee
Since $u_0 \in \MM_{\infty}$, therefore from Lemma \ref{Tab blowup are global}, $\TT_{\sqrt{a_+}, \sqrt{a_-}} (u_0) \in \PP_{\infty}$ (c.f. Definition \ref{global solution}). From the representation theorem \cite[Theorem 7.3]{ACF84}, $$|\partial \{\TT_{\sqrt{a_+}, \sqrt{a_-}}(u_0)>0\}\cap B_R^+|= |\partial \{u_0>0\}\cap B_R^+| =0.$$From \eqref{c+}, \eqref{c-} and the fact that $|\partial \{u_0>0\}\cap B_R^+| =0$, we obtain the claim \eqref{pw convergence}. Since $|\chi_{\{u_k>0\}\cap B_R^+}|\le 1$, the claim \eqref{L1 convergence} follows from Lebesgue's dominated convergence theorem.
\end{proof}

Now, let us prove Theorem \ref{tangential touch} for $u\in \MM_{1}$. That is the free boundaries of functions in $\MM_{1}$ touch the fixed boundary in tangential fashion and the modulus of continuity $\sigma$ in Theorem \ref{tangential touch} is independent of the choice of $u\in \MM_{1}$. 
In this direction, we prove that for every cone $K_{\eps}$ (as defined in \eqref{Keps}), there exists a radius $r_{\eps}$ such that $F(u)\cap B_{r_{\eps}}^+ \subset B_2^+ \setminus K_{\eps}$. 

\begin{lemma}[Tangential touch for functions in $\MM_{1}$]\label{tangential touch with density}
Let $u\in \MM_{1}$, then for all $\eps>0$ there exists $\rho(\eps, \mu, \alpha, \DD, \lambda_{\pm}, N)>0$ such that 
$$
F(u_0)\cap B_{\rho}^+ \subset  B_{\rho}^+ \setminus K_{\eps}
$$
where 
\be\label{Keps}
K_{\eps}= \big \{x\,:\, x_N \ge \eps \sqrt{x_1^2+...+x_{N-1}^2}  \big \}.
\ee
\end{lemma}

\begin{proof}
We assume by contradiction that the free boundaries do not touch the origin in tangential fashion. Then there exists $\eps>0$ and a sequences $v_j\in \MM_{1}$ and $x_j\to 0$ such that $x_j\in F(v_j)\cap K_{\eps}$ for all $j\in \N$. Let $r_j:=|x_j|$ and we define the blowups $u_j := (v_j)_{r_j}$.  

Let $u_0 = \lim_{j\to \infty}u_j$ be a limit as in Proposition \ref{blowup convergence}. Also let $x_0\in B_1^+ \cap K_{\eps}$ be a subsequential limit (for a subsequence still called $x_j$) such that $x_0= \lim_{j\to \infty}\frac{x_j}{|x_j|}$, since $v_j(x_j)=0$ therefore on rescaling, $u_j(\frac{x_j}{r_j})=\frac{1}{r_j}v_j(x_j)=0$.  Since $u_j \to u_0$ in $L^{\infty}(B_R^+)$ (c.f. Proposition \ref{blowup convergence}), we have
$$
u_0(x_0)=\lim_{j\to \infty} u_j \left (\frac{x_j}{r_j} \right )=0.
$$
We can see that $x_0\in \partial B_1\cap K_{\eps}$. Since $v_j \in \MM_{1}$, from Lemma \ref{L1 convergence} and \ref{P6}, for any given $R>0$ we have 
\be\label{density on rescaling}
\begin{split}
\frac{|\{u_0>0\}\cap B_R^+|}{|B_R^+|}= \fint_{B_R^+}\chi_{\{u_0>0\}}\,dx &=\lim_{j\to \infty}\fint_{B_R^+}\chi_{\{u_j>0\}}\,dx\\
&=  \lim_{j\to \infty} \frac{1}{|B_{R r_j}^+|} \int_{B_{Rr_j}^+}\chi_{\{v_j>0\}}\,dx\\
& =  \lim_{j\to \infty}\frac{|\{v_j>0\}\cap B_{R r_j}^+|}{|B_{R r_j}^+|} > \DD.
\end{split}
\ee
From Remark \ref{prop of Tab}, $\{u_0>0\} = \{  \TT_{\sqrt{a_+},\sqrt{a_-}} (u_0)>0\}$. Therefore we can write 
\be\label{density 1}
\frac{|\{ \TT_{\sqrt{a_+},\sqrt{a_-}}(u_0)>0\}\cap B_R^+|}{|B_R^+|}>\DD,\qquad \forall R>0.
\ee
The computations done in \eqref{density on rescaling}, in fact shows that the density property remains invariant under blowup of any function. That is why we can show that for the function $\Big [ \TT_{\sqrt{a_+},\sqrt{a_-}}(u_0)\Big ]_0$ which is the blowup of $ \TT_{\sqrt{a_+},\sqrt{a_-}}(u_0)$ we have 
\be\label{density 2}
\frac{\Big |\Big  \{\Big [  \TT_{\sqrt{a_+},\sqrt{a_-}}(u_0)\Big ]_0>0 \Big  \}\cap B_R^+\Big  |}{|B_R^+ |}>\DD,\qquad \forall R>0.
\ee
Now, we note that $u_0\in \MM_{\infty}$, therefore $ \TT_{\sqrt{a_+},\sqrt{a_-}}(u_0)$ belongs to $\PP_{\infty}$ (c.f. Lemma \ref{Tab blowup are global}). Since, from \eqref{density 1} $ \TT_{\sqrt{a_+},\sqrt{a_-}}(u_0)\not\equiv 0$, therefore from \cite[Theorem 4.2, Lemma 4.3]{KKS06} we have $ \TT_{\sqrt{a_+},\sqrt{a_-}} (u_0)\ge 0$. Moreover, from \eqref{density 2}, $\left [ \TT_{\sqrt{a_+},\sqrt{a_-}} (u_0) \right ]_0 \not \equiv 0$. Therefore, from \cite[Theorem 4.9]{KKS06}, $ \TT_{\sqrt{a_+},\sqrt{a_-}} (u_0)(x)=c\, x_N^+$ for some constant $c>0$. Thus, the function $ \TT_{\sqrt{a_+},\sqrt{a_-}}(u_0)$ cannot be equal to zero in $B_R^+$. But we have $x_0 \in \partial B_1^+ \setminus K_{\eps}$ and $ \TT_{\sqrt{a_+},\sqrt{a_-}}(u_0)(x_0)=0$. This leads us to a contradiction.

\end{proof}

 \newcommand{\noop}[1]{}

\end{document}